\documentclass[preprint,3p]{elsarticle}
\usepackage{amsmath}
\usepackage{amssymb}
\usepackage{amsthm}
\newtheorem{thm}{Theorem}[section]
\newtheorem{lem}[thm]{Lemma}
\newtheorem{cor}[thm]{Corollary}
\newtheorem{prop}{Proposition}[section]
\newdefinition{rmk}{Remark}[section]
\newproof{pf}{Proof}
\newproof{pot}{Proof of Theorem \ref{thm2}}
\newproof{poot}{Proof of Corollary \ref{co1}}
\numberwithin{equation}{section}
\newdefinition{ex}{Example}[section]
\journal{}
\theoremstyle{definition}

\newtheorem{rem}[thm]{Remark}
\begin{document}
\begin{frontmatter}

\title{  New  a priori estimates for semistable solutions of semilinear elliptic equations}
\author{A. AGHAJANI}
\ead{aghajani@iust.ac.ir}.

\address{School of Mathematics, Iran University of Science and Technology, Narmak, Tehran, Iran.}
\address{School of Mathematics, Institute for Research in Fundamental Sciences (IPM), P.O.Box: 19395-5746, Tehran, Iran.}

\begin{abstract}
We consider the  semilinear elliptic equation $-L u = f(u)$ in a  general smooth bounded domain $\Omega \subset R^{n}$
with zero Dirichlet boundary condition,  where $L$ is a uniformly elliptic operator and $f$ is a $C^{2}$ positive, nondecreasing and
convex function in $[0,\infty)$ such that $\frac{f(t)}{t}\rightarrow\infty$ as $t\rightarrow\infty$.
We prove that if $u$ is a positive semistable solution then for every $0\leq\beta<1$ we have
$$f(u)\int_{0}^{u}f(t)f''(t)~e^{2\beta\int_{0}^{t}\sqrt{\frac{f''(s)}{f(s)}}ds}~dt\in L^{1}(\Omega),$$
by a constant independent of $u$. As we shall see, a large number of  results in the literature concerning a priori bounds are immediate consequences of this estimate. In particular, among other results, we establish a priori $L^{\infty}$ bound  in dimensions $n\leq 9$, under the extra assumption that $\limsup_{t\rightarrow\infty}\frac{f(t)f''(t)}{f'(t)^{2}}< \frac{2}{9-2\sqrt{14}}\cong 1.318$. Also, we establish a priori $L^{\infty}$ bound  when $n\leq 5$ under the very weak assumption that, for some $\epsilon>0$, $\liminf_{t\rightarrow\infty}\frac{(tf(t))^{2-\epsilon}}{f'(t)}>0$ or $\liminf_{t\rightarrow\infty}\frac{t^{2}f(t)f''(t)}{f'(t)^{\frac{3}{2}+\epsilon}}>0$.
\end{abstract}
\begin{keyword} Regularity of stable  solutions; Semilinear elliptic equations; Nonlinear eigenvalue problem.
\\
\textbf{MSC(2010)}. 35K57; 35B65; 35J61

\end{keyword}

\end{frontmatter}
\section{Introduction}
This article is devoted to the study of positive semistable solutions of the following boundary value problem

\begin{equation}
\left\{\begin{array}{ll} L u+ f(u)=0& {\rm }\ x\in \Omega,\\~~u=0& {\rm }\ x\in \partial \Omega,
\end{array}\right.
\end{equation}
where $\Omega\subset R^{n}$ ($n\geq 2$) is a smooth bounded domain, $f\in C^{2}$ and $Lu:=\partial_{i}(a^{ij}(x)u_{j})$ is uniformly elliptic, namely $(a^{ij}(x))$ is a symmetric $n\times n$ matrix with bounded measurable coefficients, i.e., $a^{ij}=a^{ji}\in L^{\infty}(\Omega)$, for which there exist constants $c_{0}$ and  $C_{0}$ such that
\begin{equation}
c_{0}|\xi|^{2}\leq a^{ij}(x) \xi_{i}\xi_{j}\leq C_{0}|\xi|^{2},~~for~all~\xi\in R^{n},~x\in\Omega.
\end{equation}
By the semistability of a solution $u$ (see \cite{CSS}),   we mean that the lowest Dirichlet eigenvalue
of the linearized operator at $u$  is nonnegative. That is,
\begin{equation}
 \int_{\Omega}f'(u)\eta^{2}dx\leq \int_{\Omega} a^{ij}(x) \eta_{i}\eta_{j} dx,~~for~all~\eta\in H^{1}_{0}(\Omega).
\end{equation}
Replacing $f$ with $\lambda f$ ($\lambda\geq0$), where $f$ satisfies the assumption
\begin{equation}
f(0)>0,~~ f'\geq0 ~~\text{and} ~\lim_{s\rightarrow \infty}\frac{f(s)}{s}=\infty,
\end{equation}
then it is well known (\cite{CSS,CR,Dup}) that there exists a finite positive extremal parameter $\lambda^{*}$ such that semistable solutions exist for $\lambda\in(0,\lambda^{*})$.\\
The problem of finding a priori bounds for solutions of (1.1) under the assumption (1.4)  has been studied extensively in the literature [2-12, 15, 16] and it is shown that it depends strongly on the dimension $n$ and nonlinearity $f$.  In the  case where $L=\Delta$ and $f$ is convex, Nedev in \cite{N} obtained the $L^{\infty}$ bound for $n=2,3$ (which also holds for general $L$). When $2\leq n \leq 4$ and $L=\Delta$, the best
known result was established by Cabr\'{e} \cite{C1} who showed that the $L^{\infty}$ bound holds for arbitrary nonlinearity $f$ if in addition $\Omega$ is convex. Applying the  main estimate used
in the proof of the results  of \cite{C1}, Villegas \cite{V} got the same result replacing the condition that $\Omega$ is convex with $f$ is convex.
However, it is still an \textbf{open problem} to establish an $L^{\infty}$ estimate in dimensions $5\leq n \leq 9$, even in the case of convex domains $\Omega$ and convex nonlinearities satisfying (1.4).\\
By imposing extra assumptions on the nonlinearity $f$ much more is known, see \cite{CSS}. Let $f$ is convex and define
\begin{equation}
\tau_{-}:=\liminf_{t\rightarrow\infty} \frac{f(t)f''(t)}{f'(t)^{2}}\leq \tau_{+}:=\limsup_{t\rightarrow\infty} \frac{f(t)f''(t)}{f'(t)^{2}}.
\end{equation}
 Crandall and Rabinowitz \cite{CR} proved an a priori $L^{\infty}$ bound for semistable solutions when $0<\tau_{-}\leq \tau_{+}<2+\tau_{-}+\sqrt{\tau_{-}}$ and $n<4+2\tau_{-}+4\sqrt{\tau_{-}}$. This result was improved by  Ye and Zhou in \cite{YZ} and Sanch\'{o}n in \cite{S1} establishing that $u\in L^{\infty}$ when $\tau_{-}>0$ and $n<6+4\sqrt{\tau_{-}}$ (note that $0\leq \tau_{-}\leq1$ always hold by the assumptions on $f$). Moreover if $0<\tau_{-}\leq \tau_{+}<1$ then using an iteration argument in \cite{CR} one can show that $u\in L^{\infty}$ whenever $n<2+\frac{4}{\tau_{+}}(1+\sqrt{\tau_{-}})$. In \cite{S1} Sanch\'{o}n proved that $u\in L^{\infty}$ whenever $\tau_{-}=\tau_{+}\geq0$ and $n\leq 9$. As we have seen all the above results and others results in the literature considering $\tau_{-}$ and $\tau_{+}$ assume $\tau_{-}>0$. However, recently Cabr\'{e},  Sanch\'{o}n and Spruck \cite{CSS} proved  interesting results without assuming  $\tau_{-}>0$ and any lower bound on $f'$ nor any bound on $f''$. They considered in \cite{CSS} convex nonlinearities $f\in C^{2}$ satisfying (1.4) and one of the following conditions:\\
For every $\epsilon>0$ there exist $T_{\epsilon}$ and $C_{\epsilon}$ such that
\begin{equation}
f'(t)\leq C_{\epsilon} f(t)^{1+\epsilon}~~~ for~all~~t>T_{\epsilon},
\end{equation}
or,\\
there exist $\epsilon>0$,  $T_{\epsilon}$ and $C_{\epsilon}$ such that
\begin{equation}
f'(t)\leq C_{\epsilon} f(t)^{1-\epsilon}~~~ for~all~~t>T_{\epsilon}.
\end{equation}
They showed in \cite{CSS} that, under  condition (1.6) $u\in L^{\infty}$ when $n\leq 5$, and for $n\geq6$, $u\in W^{1,p}_{0}(\Omega)$  for all $p<\frac{n}{n-5}$. In particular, if $n\leq 9$ then $u\in H_{0}^{1}(\Omega)$.\\
Also, under condition (1.7) they showed that
$u\in L^{\infty}$ when $n< 6+\frac{4\epsilon}{1-\epsilon}$, and if $n\geq 6+\frac{4\epsilon}{1-\epsilon}$ then $u\in W^{1,p}_{0}(\Omega)$  for all $p<\frac{(1-\epsilon)n}{(1-\epsilon)n-5+3\epsilon}$. In particular, if $n\leq 10+\frac{4\epsilon}{1-\epsilon}$ then $u\in H_{0}^{1}(\Omega)$. As a corollary they proved the following results
\begin{equation}
if ~\tau_{+}<1~ and ~n<2+\frac{4}{\tau_{+}}~ then ~u\in L^{\infty}
\end{equation}
and
\begin{equation}
if ~\tau_{+}=1~ and~n<6~ then ~u\in L^{\infty}.
\end{equation}

 Note that in both the above results (also in the rest of this paper), \textbf{$u\in L^{p}(\Omega)$ or $u\in W^{1,p}(\Omega)$ mean that $u$ is bounded in $L^{p}(\Omega)$ or $ W^{1,p}(\Omega)$ by a constant independent of $u$}. Also, throughout the paper $C$  is a generic constant independent of $u$, which may take different values
in different places.\\

In this paper, we improve most of the above results by proving  the following main results using the semistability inequality (1.3) and a standard regularity result for uniformly elliptic equations.
\begin{thm}
Let $f\in C^{2}$ be convex and satisfy $(1.4)$. Let $u$ be a positive semistable solution of problem $(1.1)$. Then for every $0\leq\beta<1$ we have

\begin{equation}
H_{f,\beta}(u):=f(u)\int_{0}^{u}f(t)f''(t)~e^{2\beta\int_{0}^{t}\sqrt{\frac{f''(s)}{f(s)}}ds}~dt \in L^{1}(\Omega).
\end{equation}
\end{thm}
\begin{thm}
Let  $u\in H^{1}_{0}(\Omega)$ be a nonegative weak solution of problem $(1.1)$  with $f$ satisfies $(1.4)$. If there exists a positive constant $C$ independent of $u$ such that
\begin{equation}
||u||_{L^{1}(\Omega)}\leq C ~~and~~||\frac{\tilde{f}(u)^{\alpha}}{u^{\sigma}}||_{L^{1}(\Omega)}\leq C,~~for~some~0\leq \sigma\leq\alpha,
\end{equation}
where $\tilde{f}(u)=f(u)-f(0)$ and $\alpha\geq 1$, then
\begin{equation}
||u||_{L^{\infty}(\Omega)}\leq C~~for~~n<2\alpha.
\end{equation}
Also, if $n>2\alpha$ and $0\leq \frac{n-2}{n}\sigma< \alpha-1 $ then we have
\begin{equation}
||u||_{L^{r}(\Omega)}\leq C~~\text{for ~all} ~r<\frac{(\alpha-\sigma)n}{n-2\alpha},
\end{equation}
\begin{equation}
||f(u)||_{L^{r}(\Omega)}\leq C~~\text{for ~all} ~r<\frac{(\alpha-\sigma)n}{n-2\sigma},
\end{equation}
\begin{equation}
||u||_{W^{1,r}_{0}(\Omega)}\leq C~~ \text{for~ all} ~r<\frac{(\alpha-\sigma)n}{n-\alpha-\sigma}.
\end{equation}
In particular, if $\alpha<2+\sigma$ then
\begin{equation}
||u||_{H^{1}_{0}(\Omega)}\leq C~~ \text{for} ~n<\frac{2(\alpha+\sigma)}{2+\sigma-\alpha}.
\end{equation}
\end{thm}

Notice that, in Theorem 1.2, if $\alpha\geq2+\sigma$ then obviously we have $||u||_{H^{1}_{0}(\Omega)}\leq C$. Indeed, we then have $\frac{\tilde{f}(u)^{\alpha}}{u^{\sigma}}\geq \tilde{f}(u)^{2} \frac{\tilde{f}(u)^{\sigma}}{u^{\sigma}}$ gives $||\tilde{f}(u)^{2}||_{L^{1}(\Omega)}\leq C$ (by (1.11) and the superlinearity of $f$, i.e., $\lim_{s\rightarrow \infty}\frac{f(s)}{s}=\infty$), and as we shall see later this immediately gives $||u||_{H^{1}_{0}(\Omega)}\leq C$.

To see how the above results work and compare them with previous ones, first as an example take $f(t)=e^{t}$. Then from the estimate (1.10) we get
 $$(2+2\beta)H_{f,\beta}(u)=e^{(3+2\beta)u}-e^{u}=f(u)^{3+2\beta}-f(u)\in L^{1}(\Omega)~ for ~every ~~0\leq\beta<1,$$
 that also implies $f(u)^{3+2\beta}\in L^{1}(\Omega)$ for every $0\leq\beta<1$. Now (1.12) simply  gives  $u\in L^{\infty}(\Omega)$ when $n<10$. \\
 As an another example take $f(t)=(1+t)^{p}$, $p>1$. Then (1.10) easily gives
 \begin{equation*}
f(u)^{3-\frac{1}{p}+2\beta\sqrt{\frac{p-1}{p}}}\in L^{1}(\Omega)~ for ~every ~~0\leq\beta<1.
\end{equation*}
Now, letting $\gamma:=3-\frac{1}{p}+2\beta\sqrt{\frac{p-1}{p}}$,  then from the definition of $f$ it is easy to see that $\frac{\tilde{f}(u)^{\alpha}}{u^{\alpha}}\leq f^{\gamma}(u)\in L^{1}(\Omega)$ where $\alpha:=\frac{p}{p-1}\gamma$. Hence from (1.12) we get $u\in L^{\infty}(\Omega)$ for $n<2\alpha$, and since $\beta<1$ is arbitrary we get

$$u\in L^{\infty}(\Omega)~~for~~n<2(1+\frac{2p}{p-1}+2\sqrt{\frac{p}{p-1}}).$$
The above results are the same as results obtained by Crandall and Rabinuwitz \cite{CR}.\\ \\

Note that by the assumptions of Theorem 1.1, from (1.10) it is easy to see that $H_{f,\beta}(u)\geq C f(u)f'(u)$, that gives $f(u)f'(u)\in L^{1}(\Omega)$. This together with the fact that  $f'(u)\geq \frac{\tilde{f}(u)}{u}$ (comes from the convexity of $f$) give $\frac{\tilde{f}(u)^{2}}{u}\in L^{1}(\Omega)$. Hence, from Theorem 1.2 with $\alpha=2$ and $\sigma=1$ we get
\begin{equation}
||u||_{L^{\infty}(\Omega)}\leq C~~for~~n<4
\end{equation}
and
\begin{equation}
||u||_{H^{1}_{0}(\Omega)}\leq C~~ \text{for} ~n<6.
\end{equation}
The above results are the main results of G. Nedev in \cite{N}.

Now suppose that $\tau_{-}>0$. Then for $\tau<\tau_{-}$ there exists  $T_{\tau}$ such that $f(t)f''(t)\geq \tau f'(t)^{2}$ for $t\geq T_{\tau}$, that also gives $\frac{f''(t)}{f(t)}\geq \tau \frac{f'(t)^{2}}{f(t)^{2}}$ for $t\geq T_{\tau}$. Then using H\"{o}lder's inequality,  for $r>0$ sufficiently large we have
$$H_{f,\beta}(r)\geq C \tau f(r)\int_{0}^{r}f'(t)^{2}~e^{2\beta\sqrt{\tau}\int_{0}^{t}\frac{f'(s)}{f(s)}ds}~dt$$
$$\geq Cf(r)\int_{0}^{r}f'(t)^{2}f(t)^{2\beta\sqrt{\tau}}dt\geq Cf(r)\frac{(\int_{0}^{r}f'(t)f(t)^{\beta\sqrt{\tau}}dt)^{2}}{r}\geq C\frac{f(r)^{\delta}}{r},$$
where $\delta:=3+2\beta\sqrt{\tau}$ and $C$ is a constant independent of $u$ and depends on $\tau$ and $\beta$. Now, since $0\leq\beta<1$ and $\tau<\tau_{-}$ were arbitrary then from Theorem 1.2 we get
 \begin{equation}
u\in L^{\infty}(\Omega)~for ~~n<6+4\sqrt{\tau_{-}}.
\end{equation}
In particular, if $\tau_{-}<\frac{9}{16}$ then $u\in L^{\infty}(\Omega)$ when $n\leq 9$. Also, if $n>6+4\sqrt{\tau_{-}}$ then from (1.15) we get
\begin{equation}
||u||_{W^{1,r}_{0}(\Omega)}\leq C~~ \text{for~ all} ~r<\frac{2(1+\sqrt{\tau_{-}})n}{n-4-2\sqrt{\tau_{-}}}.
\end{equation}
The above results are the same as those obtained in \cite{S1,YZ} when $L=\Delta$.\\
Notice that, to get the above estimates (1.19-20) we assumed that $f\in C^{2}$ satisfies $\tau_{-}>0$. However, when we know only $f\in C^{1}$ then we have the following alternative.
\begin{prop}
Let  $f\in C^{1}$  satisfy $(1.4)$ and there exist  $s_{0}>0$ such that $f^{1-\delta}$ is convex in $[s_{0},\infty)$ for some $0<\delta<1$. If  $u$ is a semistable solution of  problem $(1.1)$ then
 \begin{equation}
u\in L^{\infty}(\Omega)~for ~~n<6+4\sqrt{\delta},
\end{equation}
and
\begin{equation}
||u||_{W^{1,r}_{0}(\Omega)}\leq C~~ \text{for~ all} ~r<\frac{2(1+\sqrt{\delta})n}{n-4-2\sqrt{\delta}}.
\end{equation}
\end{prop}
Note that if we assume that $f\in C^{2}$ then it is easy to see that $\tau_{-}>0$ implies that for every $0<\epsilon<1$, $f^{1-\tau_{-}+\epsilon}$ is convex in $[s_{\epsilon},\infty)$ for some $s_{\epsilon}>0$. Hence, the above result gives $u\in L^{\infty}(\Omega)$ for $n<6+4\sqrt{\tau_{-}-\epsilon}$, and since $0<\epsilon<1$ is arbitrary we get $u\in L^{\infty}(\Omega)$ for $n<6+4\sqrt{\tau_{-}}$.\\

Now consider the case  $\tau_{+}<\infty$ (we don't assume that $\tau_{-}>0$). Then the following corollary  improves the results (1.8) and (1.9).
\begin{prop}
Let  $f\in C^{2}$ be  convex and satisfy $(1.4)$, and $u$ be a positive semistable solution of  problem $(1.1)$. The following assertions hold:\\
(a) If $\tau_{+}=0$ then $u\in L^{\infty}(\Omega)$ for every $n\in \mathbb{N}$.\\
(b) If $\tau_{+}<\frac{2}{9-2\sqrt{14}}\cong 1.318$ and $n<10$ then $u\in L^{\infty}(\Omega).$\\
(c) If $0<\tau_{+}<\infty$ and $n<\max\{2+\frac{4}{\tau_{+}}+\frac{4}{\sqrt{\tau_{+}}},~4+\frac{2}{\tau_{+}}+\frac{4}{\sqrt{\tau_{+}}}\}$ then $u\in L^{\infty}(\Omega).$\\
(d) If $n>4+\frac{2}{\tau_{+}}+\frac{4}{\sqrt{\tau_{+}}}$ then
\begin{equation}
||u||_{W^{1,r}_{0}(\Omega)}\leq C~~ \text{for~ all} ~r<\frac{n}{n-3-\frac{2}{\tau_{+}}-\frac{4}{\sqrt{\tau_{+}}}}.
\end{equation}
In particular, if $\tau_{+}\neq\infty$ then $u\in H^{1}_{0}(\Omega)$ for $n<7$, and if $\tau_{+} < \frac{4}{11-4\sqrt{7}}\cong 9.592$ then $u\in H^{1}_{0}(\Omega)$ for $n<10$.
\end{prop}
Note that here we do not assume that $\tau_{+}\leq1$ as in \cite{CSS}. Also part (b) improve (1.9) from $n<6$ to $n<10$ even under the weaker assumption that $\tau_{+}<\frac{2}{9-2\sqrt{14}}$, instead of $\tau_{+}\leq1$. Indeed from part (c) we see that to get the regularity up to dimension $n<6$ we need  $\tau_{+}<10+4\sqrt{6}\simeq 19.79$. Moreover, part (c) improve (1.8) even in the case  $\tau_{+}\leq1$. Furthermore, as we have mentioned before, using previous results in the literature and an iteration argument in \cite{CR} one can show that if $0<\tau_{-}\leq \tau_{+}<1$ then   $u\in L^{\infty}$ whenever $n<2+\frac{4}{\tau_{+}}(1+\sqrt{\tau_{-}})$. However, we have
$2+\frac{4}{\tau_{+}}(1+\sqrt{\tau_{-}})\leq 2+\frac{4}{\tau_{+}}+\frac{4}{\sqrt{\tau_{+}}}$, hence part (c) of Proposition 1.2  also improve this result, without having the extra condition that $\tau_{-}>0$.\\
Also notice that from the above proposition we infer that if $\tau_{+}=\tau_{-}$ then $u\in L^{\infty}(\Omega)$ for $n\leq 9$ since in this case we must have $\tau_{+}\leq 1$ (as  $\tau_{-}\leq 1$ always holds), hence from part (c) we get $u\in L^{\infty}(\Omega)$ for $n\leq 9$.\\
If we know $0<\tau_{-}\leq \tau_{+}<\infty$ then we get better results.
\begin{cor}
Let  $f\in C^{2}$ be  convex and satisfy $(1.4)$, and $u$ be a positive semistable solution of  problem $(1.1)$. If $0<\tau_{-}\leq \tau_{+}<\infty$, then, in addition to conclusions of Proposition $1.2$ and estimates $(1.19-20)$ we also have
\begin{equation}
u\in L^{\infty}(\Omega),~~\text{for}~~n<6+\frac{4}{\sqrt{\tau_{+}}}.
\end{equation}
In particular, if $\tau_{+}<\frac{16}{9}$ then $u\in L^{\infty}(\Omega)$ for $n\leq 9$.
\end{cor}
Notice that, the above results and those in the literature including the assumption $\tau_{-}>0$,  give the uniform $L^{\infty}(\Omega)$ bound for semistable solutions at least up to dimension $6$. However, in the case when  $\tau_{-}=0$, we can use  the following consequence of Theorems 1.1 and 1.2 that gives the uniform $L^{\infty}(\Omega)$ bound  up to dimension $5$ under a very weak condition.
\begin{cor}
Let  $f\in C^{2}$ be  convex and satisfy $(1.4)$. Assume in addition that, for some $0\leq\gamma\leq 2$ and $\epsilon>0$ such that $\epsilon-\gamma>\frac{1}{2}$ we have
\begin{equation}
\liminf_{t\rightarrow\infty}\frac{t^{2-\gamma}f(t)^{1+\gamma}f''(t)}{f'(t)^{1+\epsilon}}>0.
\end{equation}
Then if $u$ is a positive semistable solution of  problem $(1.1)$, we have $||u||_{L^{\infty}(\Omega)}\leq C$ for $n\leq 5.$\\
In particular, taking $\gamma=0$ in (1.25), we see that if   for some $\epsilon>\frac{1}{2}$ we have
\begin{equation*}
\liminf_{t\rightarrow\infty}\frac{t^{2}f(t)f''(t)}{f'(t)^{1+\epsilon}}>0,
\end{equation*}
then $||u||_{L^{\infty}(\Omega)}\leq C$  for $n\leq 5.$
\end{cor}

The following preposition improves the main results of \cite{CSS} that were based on assumptions (1.6) and (1.7).
\begin{prop}
Let  $f\in C^{2}$ be  convex and satisfy $(1.4)$. Assume in addition that, for some $0\leq\gamma<\infty$ and $0\leq\delta\leq \gamma$ there exist $T:=T_{\gamma,\delta}$ and $C:=C_{\gamma,\delta}$ such that
\begin{equation}
f'(t)\leq C t^{\delta} f(t)^{\gamma}~~~ \text{for~all}~~t>T.
\end{equation}
Then if $u$ is a positive semistable solution of  problem $(1.1)$,  we have
\begin{equation}
f'(u)\in L^{1+\frac{2}{\gamma+\delta}}(\Omega)~~and~~\frac{\tilde{f}(u)^{2+\frac{1}{\gamma}}}{u^{1+\frac{1+\delta}{\gamma}}}\in L^{1}(\Omega).
\end{equation}
As a consequence
\begin{equation}
u\in L^{\infty}(\Omega)~ \text{for}~ n<\max\{ 4+\frac{2}{\gamma},~2+\frac{4}{\gamma+\delta}\}.
\end{equation}
In particular we have:\\
(i) If $\gamma<\infty$ then $u\in L^{\infty}(\Omega)$ for $n\leq 4$, and if $\gamma< 2$ then $u\in L^{\infty}(\Omega)$ for $n\leq 5$ .\\
(ii) If $n> 4+\frac{2}{\gamma}$ and $(n-2)\delta<2(\gamma+1)$ then
$$u\in L^{p}(\Omega)~ \text{for}  ~p<\frac{(1-\frac{\delta}{\gamma})n}{n-4-\frac{2}{\gamma}},~~ f(u)\in L^{p}(\Omega)~ \text{for} ~p<\frac{(1-\frac{\delta}{\gamma})n}{n-2-\frac{2(1+\delta)}{\gamma}}$$
 and
$$u\in W^{1,p}_{0}(\Omega)~~ \text{for} ~p<\frac{(1-\frac{\delta}{\gamma})n}{n-3-\frac{2+\delta}{\gamma}},$$
 $$u\in H^{1}_{0}(\Omega)~~ \text{for} ~n<\frac{6\gamma+2\delta+4}{\gamma+\delta}.$$
 In particular if $3\gamma+7\delta<4$ then $u\in H^{1}_{0}(\Omega)$ for $n\leq 9$.
\end{prop}
Note that, taking $\delta=0$ and $\gamma=1+\epsilon$ for some ($\epsilon>0$) in (1.26), then we have a weaker condition than (1.6), that we need (1.6) holds only for some $\epsilon\in(0,1)$ (not for every $\epsilon>0$ as in \cite{CSS}), but we get the regularity up to dimension $n\leq 5$. Also, from the last assertion of the above corollary we see that, if $0<\epsilon<\frac{1}{3}$ then $u\in H^{1}_{0}(\Omega)$ for $n\leq 9$. Note that, by the above corollary, to get the regularity up to dimension $n\leq 5$ we need  only to have, for some $\epsilon>0$, there exists a $T=T_{\epsilon}$ such that $$f'(t)\leq C(tf(t))^{2-\epsilon}~~\text{ for}~ t\geq T.$$
Also, taking $\delta=0$ and $\gamma=1-\epsilon$ for some ($0<\epsilon<1$) in (1.26), we have the condition (1.7). Then from (1.27-28) we get $f'(u)\in L^{\frac{3-\epsilon}{1-\epsilon}}(\Omega)$ and
$$u\in L^{\infty}(\Omega)~ \text{for}~ n<2+\frac{4}{1-\epsilon}=6+\frac{4\epsilon}{1-\epsilon}.$$
Moreover, $u\in H^{1}_{0}(\Omega)$ for $n\leq 10+\frac{4\epsilon}{1-\epsilon}$.

For example take a convex nonlinearity $f$ such that $f(t)=t\ln t$ for $t$ large. It is easy to see that $f$ satisfies (1.7) for every $0\leq \epsilon<1$, hence from (1.28) we have $u\in L^{\infty}(\Omega)$ in every dimension $n$.\\
\begin{rem}
As we have mentioned before, in dimension $n=4$, Cabr\'{e} \cite{C1} and Villegas \cite{V} showed the uniform $L^{\infty}$ bound for arbitrary nonlinearity $f$ if $\Omega$ is convex, or arbitrary domain $\Omega$ if $f$ is convex. For the proof, they used a geometric Sobolev inequality on general hypersurface of $R^{n}$ to bound the $L^{\infty}(\Omega)$ norm of every positive semistable
solution $u$ by the $W^{1,4}$ norm of $u$ on the set $\{u < t\}$  where $t$ can be
chosen arbitrarily. However, the above proposition shows that we can get the same result in dimension $n=4$ and arbitrary smooth bounded domain $\Omega$, with a more simple proof using the semistability inequality,  under the very weak extra condition that for some $\gamma<\infty$ (arbitrarily large) we have $\limsup_{t\rightarrow\infty}\frac{f'(t)}{f(t)^{\gamma}}<\infty$.
\end{rem}
Brezis and V$\acute{a}$zquez in \cite{BV}  showed that under the extra condition that
$\liminf_{t\rightarrow\infty}\frac{tf'(t)}{f(t)}>1$ or equivalently
\begin{equation}
tf'(t)-f(t)\geq\epsilon f( t),~~~t>T_{\epsilon},
\end{equation}
for some $\epsilon>0$, then
we have  $u\in H^{1}_{0}(\Omega)$. In \cite{CSS}, (1.10) is replaced with the following weaker condition that, for some $\epsilon>0$
\begin{equation}
tf'(t)-f(t)\geq\epsilon t,~~~t>T_{\epsilon}.
\end{equation}
In the following we give a weaker sufficient condition on $f$ than (1.30) to guarantee $u\in  H^{1}_{0}(\Omega)$.
\begin{prop}
Let  $f\in C^{2}$ be  convex and satisfy $(1.4)$. Assume in addition that one of the following assertions hold:\\
(i) For some $\epsilon>0$ there exists  $T=T_{\epsilon}>0$ such that
 \begin{equation}
f'(t)f\Big(t-\frac{f(t)}{f'(t)}\Big)\geq \epsilon t,~~~t>T.
\end{equation}
(ii) For some $0<\gamma<2$ there exist $C=C_{\gamma}$ and $T=T_{\gamma}$ such that
 \begin{equation}
\frac{f''(t)}{f(t)}\geq \frac{C}{t^{2}(\ln t)^{\gamma}},~~~t>T.
\end{equation}
Then if  $u$ is a positive semistable solution of  problem $(1.1)$, we have  $||u||_{H^{1}_{0}(\Omega)}\leq C$ in every dimension $n\geq 2$.
\end{prop}
Notice that from the superlinearity of $f$, i.e., $\lim_{t\rightarrow\infty}\frac{f(t)}{t}=\infty$, it is obvious  (1.31) is weaker than (1.30).
Indeed,  the left hand side of (1.31) is equal to
$\frac{f(h)}{h}(tf'(t)-f(t))$ where $h(t):=\frac{tf'(t)-f(t)}{f'(t)}\rightarrow\infty$ as $t\rightarrow\infty$ (use  L'Hospital's rule).
As an example take a nonlinearity $f$ such that  $f(t)=t (\ln t)^{a}$ for large $t$, where $0<a<1$. Then we have
$$tf'(t)-f(t)=\frac{at}{(\ln t)^{1-a}},$$
so  (1.29) or (1.30) do not hold, hence we can not apply the previous results in \cite{BV, CSS}. However, we have, for $t$ large, $f'(t)=(\ln t)^{a}+a(\ln t)^{a-1}$, hence
$$ f'(t)f\Big(t-\frac{f(t)}{f'(t)}\Big)\cong at(\ln t)^{4a-1}.$$
Thus (1.31) is satisfied if $\frac{1}{4}\leq a<1$, and by part (i) of the above proposition we have  $||u||_{H^{1}_{0}(\Omega)}\leq C$ for every domain $\Omega$ and dimension $n$. However, in this case we see that (1.32) is better than (1.31). Indeed, for $t$ sufficiently large we have
$$\frac{f''(t)}{f(t)}\cong \frac{a}{t^{2}(\ln t)^{1-2a}},$$
thus (1.32) is satisfied for every $0< a<1$. Hence, by part (ii) of the above proposition we have $||u||_{H^{1}_{0}(\Omega)}\leq C$ for every domain $\Omega$ and dimension $n$. Note that, we applied Proposition 1.4 to this example only to compare our results with previous ones, while applying Proposition 1.3 directly gives $u\in L^{\infty}(\Omega)$, implies $||u||_{H^{1}_{0}(\Omega)}\leq C$  in every dimension $n$. Indeed, here we have, for every $\gamma>0$, $f'(t)< C f(t)^{\gamma}$ for $t$ large.

\section{Preliminary estimates}
The following standard regularity result is taken from \cite{CSS}, for the proof see Theorem 3 of \cite{Ser} and Theorems 4.1 and 4.3 of \cite{Tru}, also see the explanation after Proposition 2.1 of \cite{CSS}.
\begin{prop}
Let $a^{ij}=a^{ji}$, $1\leq i,j\leq n$  be measurable functions on a bounded domain $\Omega.$ Assume that there exist positive constants $c_{0}, C_{0}$ such that $(1.2)$ holds. Let $u\in H^{1}_{0}(\Omega)$ be a weak solution of
\begin{equation}
\left\{\begin{array}{ll} L u+c(x)u=g(x)& {\rm }\ x\in \Omega,\\~~u=0& {\rm }\ x\in \partial \Omega,
\end{array}\right.
\end{equation}
with $c,g\in L^{p}(\Omega)$ for some $p\geq1$.\\
Then there exists a positive constant $ C$  independent of $u$ such that the following assertions hold:\\
(i) If $p>\frac{n}{2}$ then $||u||_{L^{\infty}(\Omega)}\leq C(|u||_{L^{1}(\Omega)}+|g||_{L^{p}(\Omega)})$.\\
(ii) Assume $c\equiv0$. If $1\leq p <\frac{n}{2}$ then $||u||_{L^{r}(\Omega)}\leq C|g||_{L^{p}(\Omega)}$ for every $1\leq r < \frac{np}{n-2p}$. Moreover, $||||_{W^{1,r}_{0}(\Omega)}\leq C$ for every $1\leq r < \frac{np}{n-p}$.
\end{prop}


The following lemma is crucial  for the proof of the main results.
\begin{lem}
Let $f\in C^{1}$  (not necessarily convex) satisfy $(1.4)$ and $g:[0,\infty]\rightarrow [0,\infty]$ be a $C^{1}$ function with $g(0)=0$ and satisfy
\begin{equation}
H(s):=g(s)^{2}f'(s)-G(s)f(s)\geq0,~~\text{for~}s~\text{sufficiently~large},
\end{equation}
where $G(s):=\int_{0}^{s}g'(t)^{2}dt$. Then if $u$ is a positive semistable solution of problem $(1.1)$, we have
$H(u)\in L^{1}(\Omega).$\\
In particular if
\begin{equation}
\limsup_{s\rightarrow\infty}\frac{G(s)f(s)}{g(s)^{2}f'(s)}<1,
\end{equation}
then
\begin{equation}
g^{2}(u)f'(u)\in L^{1}(\Omega).
\end{equation}
\end{lem}
\begin{pf}
Let $u$  be a positive semistable solution of (1.1). Take $\eta=g(u)$ as a test function in the semistability inequality  (1.3). Then we get
\begin{equation}
\int_{\Omega}a^{ij}g'(u)^{2}u_{i}u_{j}dx-\int_{\Omega} f'(u_{\lambda})g(u)^{2}dx\geq0.
\end{equation}
Now, by using the integration by part formula, we compute
\begin{equation}
\int_{\Omega}a^{ij}g'(u)^{2}u_{i}u_{j}dx=\int_{\Omega}a^{ij}u_{j}G(u)_{i}dx=-\int_{\Omega}\partial_{i}(a^{ij}u_{j})G(u)dx=\int_{\Omega} G(u)f( u)dx.
\end{equation}
Using (2.6) in (2.5) we obtain
\begin{equation}
\int_{\Omega} H(u)dx\leq0.
\end{equation}
Now from (2.2) there is an $M_{0}>0$ such that $H(s)\geq0$ for $s\geq M_{0}$, and hence using (2.7) we get
$$\int_{\Omega} |H(u)|dx=\int_{u\leq M_{0}} |H(u)|dx+\int_{u\geq M_{0}} H(u)dx\leq \int_{u\leq M_{0}} (|H(u)|-H(u))dx\leq C_{0}|\Omega|,$$
where $|\Omega|$ denotes the Lebesgue measure of $\Omega$ and $C_{0}:=\sup_{s\in[0,M_{0}]}(|H(s)|-H(s))$, and since $C_{0}$ is independent of $u$ we get $H(u)\in L^{1}(\Omega)$ that proves the first part.\\
Now suppose that (2.3)  holds and  take $\delta>0$ such that $\limsup_{s\rightarrow\infty}\frac{G(s)f(s)}{g(s)^{2}f'(s)}<\delta<1$. Then there exists an $M_{1}>0$ so that
\begin{equation}
H(s)=(1-\frac{G(s)f(s)}{g(s)f'(s)})f'(s)g(s)^{2}>(1-\delta)f'(s)g(s)^{2},~~for~~s\geq M_{1}.
\end{equation}
From (2.8) we obtain
\begin{equation*}
0\geq \int_{\Omega} H(u)dx=\int_{u<M_{1}} H(u)dx+\int_{u\geq M_{1}} H(u)dx\geq C_{1}|\Omega|+(1-\delta)\int_{u\geq M_{1}}f'(u)g(u)^{2}dx,
\end{equation*}
where  $C_{1}:=\inf_{[0,M_{1}]} H(s)$ is independent of $u$. Consequently, we have
\begin{equation*}
\int_{\Omega}f'(u)g(u)^{2}dx\leq \tilde{C}:=\Big(\frac{C_{1}}{\delta-1}+\sup_{[0,M_{1}]} f'(s)g(s)^{2}\Big)|\Omega|,
\end{equation*}
with $\tilde{C}$ independent  of $u$,  yields $g(u)^{2}f'(u)\in L^{1}(\Omega)$ that proves (2.4). \qed
\end{pf}

The following lemma will be used for the proof of Proposition 1.1.
\begin{lem}
 Let  $g:[0,\infty]\rightarrow [0,\infty]$ be a $C^{1}$ function with $g(t)>0$ for $t>0$ and there exist $s_{0}$ positive such that  $g^{\gamma}$ is convex in $[s_{0},\infty)$ for some $\gamma \in (0,1]$. If $~\limsup_{s\rightarrow\infty}\frac{f(s)}{g(s)^{2}f'(s)}=0$, then
\begin{equation}
\limsup_{s\rightarrow\infty}\frac{G(s)f(s)}{g(s)^{2}f'(s)}\leq\frac{1}{2-\gamma}\limsup_{s\rightarrow\infty}\frac{g'(s)f(s)}{g(s)f'(s)}.
\end{equation}
\end{lem}
\begin{pf}
Take $C:=\int_{0}^{s_{0}}g'(t)^{2} dt$. By the assumption $g^{\gamma}$ is convex so $g'g^{\gamma-1}$ is an increasing function in $[s_{0},\infty)$, thus for $t>t_{0}$ we can write
$$G(s)=\int_{0}^{s}g'(t)^{2} dt= C+\int_{s_{0}}^{s} [g'(t)g^{\gamma-1}(t)]~ g'(t)g^{1-\gamma}(t)dt$$
$$\leq C+ g'(s)g^{\gamma-1}(s)\int_{s_{0}}^{s} g'(t)g^{1-\gamma}(t)dt=C+\frac{g'(s)g^{\gamma-1}(s)}{2-\gamma}[g(s)^{2-\gamma}-g(s_{0})^{2-\gamma}]$$
$$\leq C+\frac{1}{2-\gamma}g(s)g'(s),$$
that easily implies (2.9).\qed
\end{pf}
\section{Proof of the main results}
{\bf Proof of Theorem 1.1}\\
Let $g$ and $H$ be as in Lemma 2.1. We write
\begin{align*}
 H(s)&=g(s)^{2}f'(s)-G(s)f(s)=f(s)\Big(g(s)^{2}\frac{f'(s)}{f(s)}-G(s)\Big):=f(s)H_{1}(s).
\end{align*}
Then from the definition of $H_{1}(s)$ we have
\begin{align}
H_{1}(s)&=\int_{0}^{s}H'_{1}(t)dt+H_{1}(0)=\int_{0}^{s}\Big(g(t)^{2}\frac{f'(t)}{f(t)}-G(t)\Big)'dt+H_{1}(0)\nonumber\\&= \int_{0}^{s}\Bigg(\frac{g(t)^{2}f''(t)}{f(t)}-f(t)^{2}\Big((\frac{g(t)}{f(t)})'\Big)^{2}\Bigg)dt+H_{1}(0).
\end{align}
Now take a $0\leq\beta<1$ and let $g(s)$ be a $C^{1}$ function with $g(0)=0$ and for some $s_{0}>0$, $g(s)=f(s)e^{\beta\int_{0}^{s}\sqrt{\frac{f''(t)}{f(t)}}dt}$ for $s>s_{0}$. Then from (3.1) we get
\begin{align}
H_{1}(s)=C+(1-\beta^{2})\int_{s_{0}}^{s}f(t)f''(t)e^{2\beta\int_{0}^{s}\sqrt{\frac{f''(t)}{f(t)}}dt},
\end{align}
and since
$$\int_{s_{0}}^{s}f(t)f''(t)e^{2\beta\int_{0}^{s}\sqrt{\frac{f''(t)}{f(t)}}dt}\geq \int_{s_{0}}^{s}f''(t)dt=f'(s)-f'(s_{0})\rightarrow\infty~as~s\rightarrow\infty,$$
from (3.2) we get
\begin{align}
H_{1}(s)\geq C\int_{0}^{s}f(t)f''(t)e^{2\beta\int_{0}^{s}\sqrt{\frac{f''(t)}{f(t)}}dt},~~\text{for~s~sufficiently~large,}
\end{align}
where $C$ is a positive constant depends only on $f$ and $\beta$. Using (3.3) and the fact that $H(u)=f(u)H_{1}(u)\in L^{1}(\Omega)$ (by Lemma 2.1) we
get
$$f(u)\int_{0}^{u}f(t)f''(t)e^{2\beta\int_{0}^{s}\sqrt{\frac{f''(t)}{f(t)}}dt}\in L^{1}(\Omega),$$
which is the desired result.
\qed\\
\begin{rem}
The following simple implication will help to simplify the proof of Theorem 1.2.
\begin{equation}
if~~f_{1}\in L^{1}(\Omega),~and~f_{2}^{q}\in L^{1}(\Omega),~~then~~(f_{1}f_{2})^{\frac{q}{q+1}}\in L^{1}(\Omega),~~(q>0).
\end{equation}
Indeed, from the assumptions we have $f_{1}^{\frac{q}{1+q}}\in L^{\frac{1+q}{q}}$ and $f_{2}^{\frac{q}{1+q}}\in L^{1+q}$, now  the H\"{o}lder inequality gives the implication.
\end{rem}
{\bf Proof of Theorem 1.2}\\
By the assumption we have $\frac{\tilde{f}(u)^{\alpha}}{u^{\sigma}}\in L^{1}(\Omega)$ for some $0\leq \sigma\leq \alpha$ ($\alpha\geq 1$), hence we have
$$\int_{\Omega} \frac{\tilde{f}(u)^{\alpha}}{u^{\alpha}}dx=\int_{u\leq1} \frac{\tilde{f}(u)^{\alpha}}{u^{\alpha}}dx+\int_{u>1} \frac{\tilde{f}(u)^{\alpha}}{u^{\sigma}}u^{\sigma-\alpha}dx\leq M|\Omega|+\int_{u>1} \frac{\tilde{f}(u)^{\alpha}}{u^{\sigma}}dx$$
$$\leq M|\Omega|+ ||\frac{\tilde{f}(u)^{\alpha}}{u^{\sigma}}||_{L^{1}(\Omega)},~~where~~M:=\sup_{0<t<1} \frac{\tilde{f}(t)^{\alpha}}{t^{\alpha}}.$$
Hence, we get $\frac{\tilde{f}(u)^{\alpha}}{u^{\alpha}}\in L^{1}(\Omega)$ or equivalently $\frac{\tilde{f}(u)}{u}\in L^{\alpha}(\Omega)$. Now similar to the proof of Corollary  2.2 in \cite{CSS} we  rewrite problem (1.1) as $Lu + c(x)u=-f(0)$ where $c(x)=\frac{\tilde{f}(u)}{u}$, hence Proposition 2.1, part (i), gives (1.12).\\
Now assume that $n>2\alpha\geq 2$ and $0\leq \frac{n-2}{n}\sigma<\alpha-1$.
From the fact that  $f(u)\in L^{1}(\Omega)$, and by the  elliptic regularity theory (see Proposition 2.1 (ii)) we get
\begin{equation}
u\in L^{q} ~~\text{for~any}~q<q_{1}:=\frac{ n}{n-2}.
\end{equation}
From (3.5) we have $(u^{\sigma})^{q}\in L^{1}$ for every $ q<\frac{q_{1}}{\sigma}$,   then using (1.11) and  Remark 3.1 we get
$$f^{ \frac{\alpha q}{1+q}}(u)=\Big(\frac{f(u)^{\alpha}}{u^{\sigma}}u^{\sigma}\Big)^{\frac{q}{1+q}}\in L^{1}(\Omega),~~for~any~q<\frac{q_{1}}{\sigma},$$
which implies that $f(u)\in L^{p}(\Omega)$ for every $1\leq p<\frac{\alpha q_{1}}{\sigma+q_{1}}$.
Note that $\frac{\alpha q_{1}}{\sigma+q_{1}}>1$ since it is equivalent to $\frac{n-2}{n}\sigma<\alpha-1$. Also we have $\frac{\alpha q_{1}}{\sigma+q_{1}}<\alpha<\frac{n}{2}.$
Again the elliptic estimates gives
\begin{equation}
u\in L^{q} ~\text{for~every}~q<q_{2}:=\frac{\alpha nq_{1}}{(\sigma+q_{1})n-2\alpha q_{1}},
\end{equation}
and by Remark 3.1 and similar as above we get $f(u)\in L^{p}$ for every $1\leq p<\frac{\alpha q_{2}}{\sigma+q_{2}}$.
Using a bootstrap procedure we can prove that $u\in L^{q}$ for every $1\leq q<q_{m}$ and $f(u)\in L^{p}$ for every $1\leq p<\frac{\alpha q_{m}}{\sigma+q_{m}}$ ($~m=1,2,...,$) where
\begin{equation}
q_{m}:= \frac{\alpha nq_{m-1}}{(\sigma+q_{m-1})n-2\alpha q_{m-1}}.
\end{equation}
Now it is easy to see that $q_{m}$ is a bounded increasing sequence with the limit $q_{\infty}$ given by
\begin{equation}
q_{\infty}=\frac{(\alpha-\sigma) n}{n-2\alpha},
\end{equation}
that proves (1.13) and (1.14). To get (1.15) it suffices to use (1.14) and Proposition 2.1, part (ii).
\qed \\
\begin{rem}
In most cases, proofs of $L^{\infty}(\Omega)$ a priori estimates in the literature are based on a uniform $L^{1}(\Omega)$ bound for functions such as $\frac{\tilde{f}(u)^{\alpha}}{u^{\sigma}}$, for some $\alpha>1$ and $0\leq\sigma\leq\alpha-1$, followed by an iterative argument of Nedev \cite{N} and standard regularity results to show $u\in L^{\infty}(\Omega)$,  for $n<2\alpha$. Our proof, however, is a direct consequence of Proposition 2.1, with an improvement of the range $0\leq\sigma\leq\alpha-1$ to $0\leq\sigma\leq\alpha$.
\end{rem}
{\bf Proof of proposition 1.1}\\
 Let $g(s)$ be a $C^{1}$ function with $g(0)=0$ and $g(s)=f^{\beta}$ for  $s>s_{0}$, where $\beta\geq 1$, and $G$ as in Lemma 2.1. Since $g^{\frac{1-\delta}{\beta}}=f(s)^{1-\delta}$ is convex, then by using (2.9) with $\gamma=\frac{1-\delta}{\beta}$ in Lemma 2.2 we have
\begin{equation}
\limsup_{t\rightarrow\infty}\frac{G(s)f(s)}{g(s)^{2}f'(s)}\leq \frac{\beta^{2}}{2\beta-1+\delta}.
\end{equation}
Now let $\beta<1+\sqrt{\delta}$, then we have  $\frac{\beta^{2}}{2\beta-1+\delta}<1$, hence form (3.9) and Lemma 2.1 we have
\begin{equation}
f(u)^{2\beta}f'(u)\in L^{1}(\Omega).
\end{equation}
From the convexity of $f$ we have $f'(t)\geq\frac{\tilde{f}(t)}{t}$, thus from (3.10) we get
\begin{equation}
\frac{\tilde{f}^{2\beta+1}}{u}\in L^{1}(\Omega),
\end{equation}
hence Theorem 1.2 gives $u\in L^{\infty}(\Omega)$ for $n<2+4\beta$. Now since $\beta<1+\sqrt{\delta}$ was arbitrary we get
 \begin{equation}
u\in L^{\infty}(\Omega)~for ~~n<6+4\sqrt{\delta},
\end{equation}
that proves (1.21). Similarly, using Theorem 2.2 and (3.11) we can prove (1.22).
\qed\\
{\bf Proof of Proposition 1.2}\\
Suppose that $\tau_{+}<\infty$. Then for $\tau>\tau_{+}$ there exists $T_{\tau}$ such that
  \begin{equation}
f(t)f''(t)\leq \tau f'(t)^{2}~~\text{for}~~ t\geq T_{\tau},
\end{equation}
 that also gives $\frac{f''(t)}{f(t)}\geq \frac{1}{\tau} \frac{f''(t)^{2}}{f'(t)^{2}}$ for $t\geq T_{\tau}$. From the convexity and superlinearity of $f$ we have $f'(t)\rightarrow \infty$ as $t\rightarrow\infty$, hence  $\int_{0}^{t}\sqrt{\frac{f''(s)}{f(s)}}ds\geq C+\frac{1}{\sqrt{\tau}} \ln f'(t)$ for large $t$. Then from the definition of $H_{f,\beta}$ in (1.10) we see that for $r>0$ sufficiently large
  \begin{equation}
H_{f,\beta}(r)\geq Cf(r)\int_{0}^{r}f(t)f''(t)f'(t)^{\frac{2\beta}{\sqrt{\tau}}}dt.
\end{equation}
Note that (3.13) is equivalent to $\frac{d}{dt}\frac{f'(t)}{f(t)^{\tau}}\leq 0$ for $t>T_{\tau}$ implies that
\begin{equation}
\frac{f'(t)}{f(t)^{\tau}}\leq C:=\frac{f'(T_{\tau})}{f(T_{\tau})^{\tau}},~~\text{for}~all~t>T_{\tau}.
\end{equation}
Using (3.15) in (3.14) we obtain, for $r$ sufficiently large
  \begin{equation}
H_{f,\beta}(r)\geq Cf(r)\int_{0}^{r}f''(t)f'(t)^{\frac{1}{\tau}+\frac{2\beta}{\sqrt{\tau}}}dt\geq Cf(r) f'(r)^{1+\frac{1}{\tau}+\frac{2\beta}{\sqrt{\tau}}},
\end{equation}
where $C$ is a constant depends on $f$ and $\beta$ but not $u$. Now if we use (3.15) in (3.16) and  Theorem 1.1 we get
\begin{equation*}
f'(u)^{1+\frac{2}{\tau}+\frac{2\beta}{\sqrt{\tau}}}\in L^{1}(\Omega),
\end{equation*}
that also gives
\begin{equation}
\frac{\tilde{f}(u)^{1+\frac{2}{\tau}+\frac{2\beta}{\sqrt{\tau}}}}{u^{1+\frac{2}{\tau}+\frac{2\beta}{\sqrt{\tau}}}}\in L^{1}(\Omega),
\end{equation}
where we used the inequality $f'(t)\geq\frac{\tilde{f}(t)}{t}$ for $t>0$. Also, using the later inequality
 in (3.16) and using Theorem 1.1 again, we get
\begin{equation}
\frac{\tilde{f}(u)^{2+\frac{1}{\tau}+\frac{2\beta}{\sqrt{\tau}}}}{u^{1+\frac{1}{\tau}+\frac{2\beta}{\sqrt{\tau}}}}\in L^{1}(\Omega).
\end{equation}
From the estimate (3.17) and Theorem 1.2 with $\alpha=\sigma=1+\frac{2}{\tau}+\frac{2\beta}{\sqrt{\tau}}$ we get
 \begin{equation}
||u||_{L^{\infty}(\Omega)}\leq C,~~for~~n<2(1+\frac{2}{\tau}+\frac{2\beta}{\sqrt{\tau}}).
\end{equation}
Now if $\tau_{+}=0$ then since  (3.19) holds for every $\tau>\tau_{+}=0$ we get $||u||_{L^{\infty}(\Omega)}\leq C$ for every $n\in \mathbb{N}$ that proves part (a). Also, if $\tau_{+}>0$ since $0\leq\beta<1$ and $\tau>\tau_{+}$ are arbitrary in (3.19) then  we get
 \begin{equation}
||u||_{L^{\infty}(\Omega)}\leq C,~~for~~n<2+\frac{4}{\tau_{+}}+\frac{4}{\sqrt{\tau_{+}}}).
\end{equation}
Also, from the estimate (3.18) and Theorem 1.2 with $\alpha=\sigma+1=2+\frac{1}{\tau}+\frac{2\beta}{\sqrt{\tau}}$ and similar as above we get
 \begin{equation}
||u||_{L^{\infty}(\Omega)}\leq C,~~for~~n<4+\frac{2}{\tau_{+}}+\frac{4}{\sqrt{\tau_{+}}}.
\end{equation}
Now, (3.20) and (3.21) prove part (c). To prove part (b), it suffices to note that for $\tau_{+}<\frac{2}{9-2\sqrt{14}}$ we have $4+\frac{2}{\tau_{+}}+\frac{4}{\sqrt{\tau_{+}}}>9$ and use part (c). Also,  using the estimate (3.18) and Theorem 2.2 we can easily prove part (d).
\qed\\
{\bf Proof of Corollaries 1.3 and 1.4}\\
Suppose $\tau_{-}>0$ then there exist a $T>0$ such that
  \begin{equation}
f(t)f''(t)\geq \frac{\tau_{-}}{2} f'(t)^{2}~~\text{for}~~ t\geq T.
\end{equation}
Now, for $\tau>\tau_{+}$ as in the proof of Proposition 1.2 and using (3.22),  (3.14) and H\"{o}lder inequality, for $r>0$ sufficiently large we have
  \begin{equation}
H_{f,\beta}(r)\geq Cf(r)\int_{0}^{r}f'(t)^{2+\frac{2\beta}{\sqrt{\tau}}}dt\geq Cf(r)\frac{(\int_{0}^{r}f'(t)dt)^{2+\frac{2\beta}{\sqrt{\tau}}}}{r^{1+\frac{2\beta}{\sqrt{\tau}}}}
\geq C\frac{f(r)^{3+\frac{2\beta}{\sqrt{\tau}}}}{r^{1+\frac{2\beta}{\sqrt{\tau}}}}.
\end{equation}
where $C$ is a constant depends on $f$ and $\beta$ but not $u$. Now similar to the proof of Proposition 1.2 and using Theorems 1.1 and 1.2 we get the desired result of Corollary 1.3.\\
To prove Corollary 1.4,  from (1.25) we deduce, there exist  $C$ and $T>0$ such that
  \begin{equation}
f(t)f''(t)\geq C \frac{f'(t)^{1+\epsilon}}{t^{2-\gamma}f(t)^{\gamma}}~~\text{for}~~ t\geq T.
\end{equation}
Hence, using (3.24) and H\"{o}lder inequality, for $r>0$ sufficiently large we have
$$H_{f,\beta}(r)\geq Cf(r)\int_{0}^{r}f(t)f''(t)dt\geq Cf(r)\int_{0}^{r}\frac{f'(t)^{1+\epsilon}}{t^{2-\gamma}f(t)^{\gamma}}\geq C\frac{f(r)}{r^{2-\gamma}f(r)^{\gamma}}\int_{0}^{r}f'(t)^{1+\epsilon}dt$$
$$\geq C\frac{f(r)}{r^{2-\gamma}f(r)^{\gamma}}\frac{(\int_{0}^{r}f'(t)dt)^{1+\epsilon}}{r^{\epsilon}}\geq C\frac{f(r)^{2+\epsilon-\gamma}}{r^{2+\epsilon-\gamma}}.$$
Hence, $\frac{\tilde{f}(u)^{2+\epsilon-\gamma}}{u^{2+\epsilon-\gamma}}\in L^{1}(\Omega)$, thus Theorem 1.2 gives $u\in L^{\infty}(\Omega)$ for $n\leq 4+2(\epsilon-\gamma)$. Now, by the assumption that $\epsilon-\gamma>\frac{1}{2}$ we get $u\in L^{\infty}(\Omega)$ for $n\leq 5$.
\qed\\
{\bf Proof of Proposition 1.3}\\
By using the assumption (1.26), for $r>0$ sufficiently large we have
$$H_{f,\beta}(r)\geq Cf(r)\int_{0}^{r}f(t)f''(t)dt\geq Cf(r)\int_{0}^{r}\frac{f'(t)^{\frac{1}{\gamma}}}{t^{\frac{\delta}{\gamma}}}f''(t)dt$$
$$\geq C\frac{f(r)}{r^{\frac{\delta}{\gamma}}}\int_{0}^{r}f'(t)^{\frac{1}{\gamma}} f''(t)dt \geq C\frac{f(r)}{r^{\frac{\delta}{\gamma}}}f'(r)^{1+\frac{1}{\gamma}},$$
that gives (using Theorem 1.2)
\begin{equation}
\frac{\tilde{f}(u)}{u^{\frac{\delta}{\gamma}}}f'(u)^{1+\frac{1}{\gamma}}\in L^{1}(\Omega).
\end{equation}
By using the inequality $f'(t)\geq \frac{\tilde{f}(t)}{t}$ for $t>0$ in (3.25) we get
\begin{equation}
\frac{\tilde{f}(u)^{2+\frac{1}{\gamma}}}{u^{1+\frac{1+\delta}{\gamma}}}\in L^{1}(\Omega).
\end{equation}
Also, from (1.26) and the superlinearity of $f$ we get ,  for $t$ large enough,  $f'(t)\leq C f(t)^{\gamma+\delta}$, that gives
$$\frac{f(t)}{t^{\frac{\delta}{\gamma}}}f'(t)^{1+\frac{1}{\gamma}}=\frac{f(t)^{\frac{\delta}{\gamma}}}{t^{\frac{\delta}{\gamma}}}
f(t)^{1-\frac{\delta}{\gamma}}f'(t)^{1+\frac{1}{\gamma}}\geq f(t)^{1-\frac{\delta}{\gamma}}f'(t)^{1+\frac{1}{\gamma}}\geq f'(t)^{1+\frac{2}{\gamma+\delta}},$$
for $t$ sufficiently large. Using the above inequality and (3.25) we deduce
\begin{equation}
f'(u)^{1+\frac{2}{\gamma+\delta}}\in L^{1}(\Omega).
\end{equation}
Now (3.26) and (3.27) prove (1.27), and an application of Theorem 1.2 completes the proof.
\qed\\
{\bf Proof of Proposition 1.4}\\
First notice that, from \cite{CSS} ( or \cite{BV} for the case $L=\Delta$), to prove the conclusion of theorem it is sufficient to show that
\begin{equation*}
uf(u)\in L^{1}(\Omega).
\end{equation*}
From the estimate (1.12) in Theorem 1.2 we have
\begin{equation}
h(u):=f(u)\int^{u}_{0}f(t)f''(t)dt\in L^{1}(\Omega).
\end{equation}
From the convexity of $f$ and Jensen's inequality we have
$$h(s)=f(s)\int^{s}_{0}f(t)d(f'(t))\geq f(s) (f'(s)-f'(0))f\Big(\frac{1}{f'(s)-f'(0)}\int^{s}_{0}td(f'(t))\Big)$$
$$=f(s) (f'(s)-f'(0))f\Big(\frac{sf'(s)-f(s)+f(0)}{f'(s)-f'(0)}\Big)\geq Cf(s) f'(s)f\Big(\frac{sf'(s)-f(s)}{f'(s)}\Big),$$
for $s$ sufficiently large. Note that by the L'Hospital's rule, we have $\lim_{t\rightarrow\infty}\frac{sf'(s)-f(s)}{f'(s)}=\infty$. Now suppose, for some $\epsilon>0$, (1.31) holds. Then from the above estimate and (3.28) we get $uf(u)\in L^{1}(\Omega))$ that proves part (i).\\
To prove part (ii), first note that from (1.32) we get, for $t$ large enough
$$\int_{0}^{t}\sqrt{\frac{f''(s)}{f(s)}}ds\geq \tilde{C} (\ln t)^{1-\frac{\gamma}{2}}.$$
By using the above inequality and (1.32)  we have, for $r$ sufficiently large
$$H_{f,\beta}(r)\geq C f(r)\int_{2}^{r}\frac{f(t)^{2}}{t^{2}(\ln t)^{\gamma}}e^{\tilde{C} (\ln t)^{1-\frac{\gamma}{2}}}dt.$$
Now using the fact that $e^{\tilde{C} (\ln t)^{1-\frac{\gamma}{2}}}\geq (\ln t)^{\gamma}$ for $t$ large enough,  the above inequality implies
\begin{equation}
H_{f,\beta}(r)\geq C f(r)\int_{2}^{r}\frac{f(t)^{2}}{t^{2}}dt\geq C \frac{f(r)}{r^{2}}\int_{0}^{r}f(t)^{2}dt.
\end{equation}
From the superlinearity of $f$ we have $\int_{0}^{r}f(t)^{2}dt\geq Cr^{3}$ for $r$ sufficiently large, hence from (3.29) we get $H_{f,\beta}(r)\geq C r f(r)$. Now, theorem 1.2 implies $uf(u)\in L^{1}(\Omega)$ that gives the desired result.
\qed \\
\section{Acknowledgement}
This research was in part supported by a grant from IPM (No. 93340123).

\end{document}